\documentclass{birkjour}
 \theoremstyle{definition}
 
 \theoremstyle{remark}

 \numberwithin{equation}{section}

\usepackage{amsfonts,amsmath,amssymb,indentfirst}
\newcommand{\er}{\mathbb{R}}

\newtheorem{Teorema}{Theorem}[section]
\newtheorem{Propriedade}[Teorema]{Proposition}
\newtheorem{Rem}[Teorema]{Remark}
\newtheorem{Lemma}[Teorema]{Lemma}

\begin{document}

%
%
%
%
%
%
%
%
%

\title[Traveling-wave solutions to a Boussinesq system]{A note on the existence of traveling-wave\\ solutions to a Boussinesq system}

\author[Filipe Oliveira]{Filipe Oliveira}

\address{Department of Mathematics, FCT - U. Nova de Lisboa, 2829-516 Caparica, Portugal}
\address{Centro de Matem\'atica e Aplica\c{c}\~oes Fundamentais, U. de Lisboa, Portugal}
\email{fso@fct.unl.pt}

\thanks{The author was partially supported by FCT (Portuguese Foundation for Science and
Technology) through the grant PEst-OE/MAT/UI0209/2011.}

\subjclass{Primary 35Q35; Secondary 76B25}

\keywords{Traveling-waves, Dispersive equations, Boussinesq systems}

\date{August 3, 2014}

\begin{abstract}
We obtain a one-parameter family $$(u_{\mu}(x,t),\eta_{\mu}(x,t))_{\mu\geq \mu_0}=(\phi_{\mu}(x-\omega_{\mu} t),\psi_{\mu}(x-\omega_{\mu} t))_{\mu\geq \mu_0}$$ of traveling-wave solutions to the Boussinesq system
\begin{displaymath}
\left\{\begin{array}{llll}
u_t+\eta_x+uu_x+c\eta_{xxx}=0\qquad (x,t)\in\er^2 \\
\eta_t+u_x+(\eta u)_x+au_{xxx}=0
        \end{array}\right.
\end{displaymath}
in the case $a,c<0$, with non-null speeds $\omega_{\mu}$ arbitrarely close to $0$ ($\omega_{\mu}\xrightarrow[\mu\to+\infty]{} 0$). We show that the $L^2$-size of such traveling-waves satisfies the uniform (in $\mu$) estimate
$\|\phi_{\mu}\|_2^2+\|\psi_{\mu}\|_2^2\leq C\sqrt{|a|+|c|},$
where $C$ is a positive constant. Furthermore, $\phi_{\mu}$ and $-\psi_{\mu}$ are smooth, non-negative, radially decreasing functions which decay exponentially at infinity. 
\end{abstract}

\maketitle
\section{Introduction}
\noindent
In \cite{Bona}, with the purpose of describing the dynamics of small-amplitude long waves propagating on the surface of an ideal fluid in a channel of constant depht, the authors introduced the four-parameter family of Boussinesq systems
\begin{equation}
 \label{Boussinesqsys}
\left\{\begin{array}{llll}
u_t+\eta_x+uu_x+c\eta_{xxx}-du_{xxt}=0 \\
\eta_t+u_x+(\eta u)_x+au_{xxx}-b_{xxt}=0.
        \end{array}\right.
\end{equation}
These systems are first-order approximations to the Euler equations in the small parameters $\alpha=\frac Ah<<1$ and $\beta=\frac{h^2}{l^2}<<1$, where $h$ is the depht of the channel (subsequently scaled to $1$), and $A$ and $\lambda$ represent
 a typical wave amplitude and a typical wavelenght respectively. Here, $\eta(x,t)$ denotes the deviation of free surface with respect to the undisturbed state (i.e. $1+\eta(x,t)$ is the total depht of the liquid at time $t$ 
and position $x$) and $u(x,t)$ is the horizontal 
velocity at height $\theta$, $0\leq \theta \leq 1$. The four parameters $a,b,c$ and $d$ are given by
\begin{equation}
 \label{abcd}
\begin{array}{llll}
\displaystyle a=\frac 12\left(\theta^2-\frac 13\right)\lambda,&\displaystyle b=\frac 12\left(\theta^2-\frac 13\right)(1-\lambda)\\
\\
\displaystyle c=\frac 12\left(1-\theta^2\right)\mu,&\displaystyle d=\frac 12\left(1-\theta^2\right)(1-\mu),
        \end{array}
\end{equation}
where, as stated in \cite{Bona}, $\lambda$ and $\mu$ are modeling parameters that do not possess a direct physical interpretation.\\ 
In \cite{Daripa}, a correction of the $c$ parameter  is proposed in order to include the contribution of the surface tension:
\begin{equation}
\label{novoc}
 c=\frac 12\left(1-\theta^2\right)\mu-\tau.
\end{equation}
The Bond number $\tau$ is given by $\displaystyle\tau=\frac{\Gamma}{\rho gh^2}$, where $\Gamma$ is the surface tension coefficient and $\rho$ the density of water.\\ 
In \cite{Chen1}, the authors prove the existence (and orbital stability) of traveling-waves to system (\ref{Boussinesqsys}) of the form
\begin{equation}
\label{traveling}
u(x,t)=\phi(x-\omega t),\eta(x,t)=\psi(x-\omega t)
\end{equation}
in the case $a,c<0$, $b=d>0$ and $ac>b^2$. Furthermore, in \cite{Chen2}, the existence of such traveling-waves with small propagation speed is obtained in the case $a,c<0$ and $b=d$.

\medskip

\noindent
In the present paper we exhibit a new family of one-parameter traveling-waves in the case $a,c<0$ and $b=d=0$. Our method has the advantage of providing radially decreasing functions and a uniform bound for the $L^2-$size of the solution. 
More precisely, we prove the following result:
\begin{Teorema}
\label{Teorema}
Let $a,c<0$. There exists a constant $\mu_0=\mu_0(a,c)$ and a one-parameter family of nontrivial traveling-wave solutions to the Boussinesq system
\begin{equation}
\label{Boussinesq}
\left\{\begin{array}{llll}
u_t+\eta_x+uu_x+c\eta_{xxx}=0 \\
\eta_t+u_x+(\eta u)_x+au_{xxx}=0
        \end{array}\right.
\end{equation}
of the form 
$$(u,\eta)=(\phi_{\mu}(x-\omega_{\mu}t), \psi_{\mu}(x-\omega_{\mu}t)),\qquad \mu\geq \mu_0,$$
with $(\phi_{\mu},\psi_{\mu})\in H^{\infty}(\er)\times H^{\infty}(\er)$, $\phi_{\mu}$ and $-\psi_{\mu}$ non-negative and radially decreasing, with exponential decay at infinity. \\
Furthermore, the speed $\omega_{\mu}$ satisfies the estimate
$$0>\omega_{\mu}>-\frac 1{C_1}\sqrt[3]{|a|+|c|}\mu^{-\frac 23}$$
and the following uniform control of the $L^2$-norm of the traveling-wave holds:
\begin{equation}
\label{l2norm}
\|\phi_{\mu}\|_2^2+\|\psi_{\mu}\|_2^2\leq C\sqrt{|a|+|c|},
\end{equation}
where $C$ is a positive constant independent of $a$, $c$ and $\mu$.
\end{Teorema}
\begin{Rem}
We stated the existence of traveling-waves $$(\phi(x-\omega t),\psi(x-\omega t)),$$ with $\phi\geq 0$ and $\psi\leq 0$, propagating with negative speed $\omega$. Noticing that $(-\phi(x+\omega t), \psi(x+\omega t))$
is also a solution to (\ref{Boussinesq}), it is straightforward to deduce the existence of non-positive traveling-waves propagating with positive speed.  
\end{Rem}
\begin{Rem}
As stated above, in \cite{Chen2}, the authors also establish the existence of traveling-waves with small speed in the particular case $a,c<0$ and $b=d=0$, although it is not clear, with the method used, if the solutions have a sign,
are radially-decreasing or if an uniform (in the speed $\omega$) estimate such as (\ref{l2norm}) holds. On the other hand, of course, the method used in \cite{Chen2} has the important advantage of covering the case $b=d\neq 0$. Either way, it
does not seem obvious to prove or to disprove that the traveling-waves found in both papers are the same.
\end{Rem}
\begin{Rem}
Note, in view of (\ref{abcd}) and (\ref{novoc}), that the case treated in Theorem \ref{Teorema} corresponds to $\lambda=\mu=1$, that is
$$a=\frac 12\left(\theta^2-\frac 13\right)\quad\textrm{ and }\quad c=\frac 12(1-\theta^2)-\tau.$$
For $\theta^2\to {\frac 13}$ and $\tau\to {\frac 13}$ we get $a,c\to 0$. The estimate $(\ref{l2norm})$ then suggests that the traveling-wave solutions vanish in this regime.\\
This is consistent with the known fact (see \cite{Daripa}) that in the case $\frac 13-\tau=\mathcal{O}(\beta)$, $\beta\to 0$ it is necessary to introduce higher order terms in the Boussinesq approximation of the Euler 
equations in order to model solitary waves.
\end{Rem}
\section{The minimization problem}
\noindent
In what follows we fix $a,c<0$. In order to prove Theorem \ref{Teorema} we introduce a variational problem whose minimizers, up to a rescaling,  correspond to a traveling-wave profile. Let us first note that 
$$(u,\eta)=(\phi(x-\omega t), \psi(x-\omega t)),\quad \phi(y),\psi(y)\xrightarrow[y\to\infty]{}0$$
is a solution to the Boussinesq system (\ref{Boussinesq}) if and only if $\phi$ and $\psi$ satisfy the stationary equation
\begin{equation}
\label{estacionaria}
\left\{
\begin{array}{llll}
\displaystyle a\phi''+\phi-\omega\psi+\phi\psi=0\\
\displaystyle c\psi''+\psi-\omega\phi+\frac 12 \phi^2=0.
\end{array}\right.
\end{equation}
For $\mu>0$, we set
$$X_{\mu}=\{(f,g)\in H^1(\er)\times H^1(\er)\,:\,N(f,g)=\|f\|_2^2+\|g\|_2^2=\mu\}$$
and consider the minimization problem
$m(\mu)=\inf \{\tau(f,g)\,:\,(f,g)\in X_{\mu}\},$
where $$\tau(f,g)=-a\int {f'}^2-c\int {g'}^2+\int f^2g+2\int fg.$$
\begin{Propriedade}
\label{ninfinito}
For all $\mu>0$, $m(\mu)>-\infty$.\\
More precisely, $\displaystyle m(\mu)\geq -\frac{C}{\sqrt[3]{|a|}}\mu^{\frac{10}3}-2\mu$, where $C$ is a positive constant.
\end{Propriedade}
\noindent
{\bf Proof}
Let $(f,g)\in X_{\mu}$. One only has to notice that by the Cauchy-Schwarz and Gagliardo-Nirenberg inequalities,
$$\left|\int fg\right|\leq \|f\|_2\|g\|_2\leq \mu
\quad\textrm{and}$$
$$\quad\left|\int f^2g\right|\leq \|f\|_4^2\|g\|_2\leq \|f'\|_2^{\frac 12}\|f\|_2^{\frac 32}\|g\|_2\leq \mu^{\frac 52}\|f'\|_2^{\frac 12}.$$
Hence,
$$\tau(f,g)\geq -a\|f'\|_2^2-\mu^{\frac 52}\|f'\|_2^{\frac 12}-2\mu=P(\|f'\|_2),$$
which is enough to conclude since $\frac{1}{\sqrt[3]{|a|}}\left(16^{-\frac 13}-4^{-\frac 13}\right)\mu^{\frac {10}3}-2\mu$ is the minimum of $P(x)=-ax^4-\mu^{\frac 52}x-2\mu$.\hfill $\blacksquare$
\begin{Propriedade}
\label{estimativaboa}
For all $\mu>0$, $\displaystyle m(\mu)<-\frac{C}{\sqrt[3]{|a|+|c|}}\mu^{\frac 53}$, where $C$ is a positive constant.
\end{Propriedade}
\noindent
{\bf Proof}
We fix a non-negative function $h\in H^1(\er)$ such that $\|h\|_2=1$ and we put $h_{\lambda}(x)=\lambda h(\lambda^2x)$, where $\lambda$ will be chosen later. 
Then, for all $\mu>0$, $\displaystyle (f,g)=\frac 1{\sqrt{2}} \left(  \mu^{\frac 23}h_{\lambda}(\mu^{\frac 13}x),-\mu^{\frac 23}h_{\lambda}(\mu^{\frac 13}x)  \right)\in X_{\mu}$ and 
$$m(\mu)\leq \tau(f,g)=\frac{|a|+|c|}{2}\mu^{\frac 53}\int (h_{\lambda}')^2-\frac 1{2\sqrt{2}}\mu^{\frac 53}\int h_{\lambda}^3-\mu\int h_{\lambda}^2$$
$$\leq \mu^{\frac 53}\left(\frac{|a|+|c|}{2}\|h_{\lambda}'\|_2^2-\frac 1{2\sqrt{2}}\|h_{\lambda}\|_3^3 \right)$$
$$\leq \mu^{\frac 53}\lambda\left(\lambda^3(|a|+|c|)\|h'\|_2^2-\frac 1{2\sqrt{2}}\|h\|_3^3\right).$$
We conclude the proof by choosing $\lambda=\frac{\epsilon}{\sqrt[3]{|a|+|c|}}$, with $\epsilon$ such that\\ $-C=\lambda^3(|a|+|c|)\|h'\|_2^2-\frac 1{2\sqrt{2}}\|h\|_3^3<0$.\hfill $\blacksquare$

\section{Existence of Minimizers}
\noindent
Let $\mu>0$ and $(f_n,g_n)$ a minimizing sequence for $m(\mu)$. By denoting $f^*$ the Schwarz symmetrization of $|f|$, it is well known that 
\begin{displaymath}
\|{f^*}'\|_2\leq \|f'\|_2,\quad \|f^*\|_2=\|f\|_2,\quad \int fg\leq \int f^*g^*\textrm{ and }\int {f}^2g\leq \int {f^*}^2g^*.
\end{displaymath}
Hence, $\tau(f^*,-g^*)\leq \tau(f,g)$ and $(f,g)\in X_{\mu}$ implies that $(f^*,-g^*)\in X_{\mu}$. Therefore we can choose a minimizing sequence $(f_n,g_n)$ with $f_n\geq$, $g_n\leq 0$ and $f_n$, $-g_n$ radially decreasing.

\medskip

\noindent
We will now apply the concentration-compactness method (\cite{Lions1},\cite{Lions2}) to prove the compacity, up to translations, of the sequence $(f_n,g_n)$ in $L^2(\er)\times L^2(\er)$ strong.\\
Following this method, we set the concentration function $\rho_n=f_n^2+g_n^2$ and put $\displaystyle Q_n(t)=\sup_{y\in\er}\int_{y-t}^{y+t}\rho_n$.
Also, we set $\displaystyle Q(t)=\lim_{n\to+\infty}Q_n(t)$  and $\displaystyle \Omega=\lim_{t\to+\infty} Q(t)$.

\medskip

\noindent
We start by ruling out vanishing:
\begin{Propriedade}
\label{previous}
 There exists $\mu_0=\mu_0(a,c)$ such that for $\mu\geq \mu_0$, $\Omega>0$.
\end{Propriedade}
\noindent
{\bf Proof}
Assume that $\Omega=0$. Since $Q(t)$ is non-negative and non-increasing, for all $t$, $Q(t)=0$. Hence, 
$$\lim_{n\to +\infty}\sup_{y\in\er}\int_{y-t}^{y+t}f_n^2=\lim_{n\to +\infty}\sup_{y\in\er}\int_{y-t}^{y+t}g_n^2=0.$$
From the Proof of \ref{ninfinito} one can infer that $(f_n)$ (and $g_n$) is bounded in $H^1(\er)$. Arguing as in \cite{Lions1} (Lemma I.1), $\|f_n\|_4\to 0$, hence $\int f_n^2g_n\to 0$ by Cauchy-Schwarz. Furthermore,
$$m(\mu)=\lim_{n\to +\infty} \left(-a\int {f_n'}^2-c\int {g_n'}^2+\int f_n^2g_n+2\int f_ng_n\right)$$   
$$\geq 2\lim_{n\to +\infty}\int f_ng_n\geq -2\mu,$$
which contradicts Proposition \ref{estimativaboa} for $\mu\geq \mu_0=\frac{2^{\frac 32}}{C_1^{\frac 32}}\sqrt{|a|+|c|}$.\hfill $\blacksquare$

\medskip

\noindent
Next, we rule out dichotomy, that is $\displaystyle 0<\Omega<\lim_{n\to+\infty}\int \rho_n.$ It is sufficient to prove the following lemma:
\begin{Lemma}
\label{theta}
 For all $\mu\geq \mu_0$ and for all $\theta>1$, $m(\theta\mu)<\theta m(\mu)$.
\end{Lemma}
\noindent
{\bf Proof}
We have
$\displaystyle \tau(\theta^{\frac 12}f_n,\theta^{\frac 12}g_n)=\theta \tau(f_n,g_n)-(\theta^{\frac 32}-\theta)\int |f_n|^2|g_n|.$\\
Also, there exists $\delta>0$ such that for all $n$ large enough, $\int |f_n|^2|g_n|\geq \delta$. Otherwise, up to a subsequence, $\int f_n^2g_n\to 0$, which is absurd for $\mu\geq \mu_0$, as seen in the previous proof of 
Proposition \ref{previous}. Finally, 
\begin{displaymath}
\displaystyle m(\theta\mu)\leq\tau(\theta^{\frac 12}{f_n},\theta^{\frac 12}{g_n})\leq\theta \tau(f_n,g_n)-\delta(\theta^{\frac 32}-\theta),
\end{displaymath}
which yields the result:\\
$\displaystyle m(\theta\mu)\leq\lim_{n\to+\infty} \theta \tau(f_n,g_n)-\delta(\theta^{\frac 32}-\theta)=\theta m(\mu)-\delta(\theta^{\frac 32}-\theta)<\theta m(\mu).$ \hfill$\blacksquare$

\medskip

\noindent
It is standard, from Lemma \ref{theta}, to prove the strict subadditivity of $m$, that is $$\forall \mu\geq \Omega,\quad m(\mu)<m(\Omega)+m(\mu-\Omega),$$ 
(see for instance Lemma 2.3 in \cite{Ohta}) which is well-known to rule out dichotomy.

\medskip

\noindent
Hence, by Lions' Theorem, we are in the compactness situation. There exists a sequence $(y_n)$ such that, up to a subsequence, $(f_n(.-y_n),g_n(.-y_n))$ converges strongly in $L^2(\er)\times L^2(\er)$ to some $(\tilde{f},\tilde{g})\in X_{\mu}$.\\
Since $(f_n,g_n)$ is bounded in $H^1(\er)$, using the compact embedding $$H^1_{rad}(\er)\hookrightarrow L^4(\er),$$ up to a subsequence, $f_n$ (respectively $-g_n$) converges strongly in $L^4$ to some radial non-negative function $f$ 
(respectively $-g$).\\
Furthermore, $(f_n,g_n)\rightharpoonup(f,g)$ in $H^1(\er)\times H^1(\er)$ weak.\\
Since 
$$\|f_n\|_2^2+\|g_n\|_2^2=\|f_n(\cdot-y_n)\|_2^2+\|g_n(\cdot-y_n)\|_2^2\xrightarrow[n\to +\infty]{}\|\tilde{f}\|_2^2+\|\tilde{g}\|_2^2=\mu,$$
we have in fact that $(f_n,g_n)\to (f,g)$ in $L^2(\er)\times L^2(\er)$ strong and, in particular, $(f,g)\in X_{\mu}$. Hence, 
$$\int f_ng_n\to \int f^2g,\,\int f_ng_n \to \int fg \textrm{ and }$$
$$\lim_{n\to +\infty} |a|\int {f_n'}^2+|c|\int {g_n'}^2 \leq |a|\int {f'}^2+|c|\int {g'}^2.$$
From these inequalities, we deduce that $\displaystyle \tau(f,g)\leq \lim_{n\to+\infty} \tau(f_n,g_n)=m(\mu).$
Since $(f,g)\in X_{\mu}$, we have in fact $\tau(f,g)=m(\mu)$ and $(f,g)$ is a minimizer for $m(\mu)$.

\section{End of the Proof of Theorem \ref{Teorema}}
\noindent
There exists a Lagrange multiplier $\lambda \in \er$ such that $\nabla \tau(f,g)=\lambda \nabla N(f,g)$, that is
\begin{equation}
 \label{min1}
\left\{
\begin{array}{llll}
\displaystyle af''+fg+g&=&\lambda f\\
\displaystyle cg''+\frac 12f^2+f&=&\lambda g. 
\end{array}\right.
\end{equation}
Multiplying these equations by $f$ and $g$ respectively and integrating by parts leads to 
\begin{displaymath}
\left\{
\begin{array}{llll}
\displaystyle -a\int(f')^2+\int f^2g+\int fg&=&\lambda \int f^2\\
\\
\displaystyle -c\int(g')^2+\frac 12\int f^2g+\int fg&=&\lambda \int g^2, 
\end{array}\right.
\end{displaymath}
and, adding the equalities,
\begin{equation}
 \label{estimativaL}
-\lambda \mu=-\tau(f,g)-\frac 12 \int f^2g\geq -\tau(f,g)=-m(\mu). 
\end{equation}
Note that, in particular, $\lambda<0$. Setting 
$$\phi(x)=-\frac 1{\lambda}f\left(\frac x{\sqrt{-\lambda}}\right)\quad\textrm{ and }\psi(x)=-\frac 1{\lambda}g\left(\frac x{\sqrt{-\lambda}}\right),$$
we obtain
\begin{equation}
 \label{min2}
\left\{
\begin{array}{llll}
\displaystyle a\phi''+\phi+\phi\psi-\frac 1{\lambda}\psi&=&0\\
\displaystyle c\psi''+\psi+\frac 12\phi^2-\frac 1{\lambda}\phi&=&0,
\end{array}\right.
\end{equation}
that is, $(\phi,\psi)$ is a solution to (\ref{estacionaria}) with speed $\displaystyle \omega=\frac 1{\lambda}.$

\medskip
\noindent
To obtain the $L^2$ size of this solution, a simple computation shows that
$$\|\phi\|_2^2+\|\psi\|_2^2=\frac{\mu}{|\lambda|^{\frac 32}}.$$
Also, from (\ref{estimativaL}) and Proposition \ref{estimativaboa}, $|\lambda| \mu \geq -m(\mu)\geq \frac{C}{\sqrt[3]{|a|+|c|}}\mu^{\frac 53},$
from where we conclude that
\begin{displaymath}
\|\phi\|_2^2+\|\psi\|_2^2\leq C\sqrt{|a|+|c|},
\end{displaymath}
where $C$ is yet another positive constant independent of $a$, $c$ and $\mu$. Also, note that 
\begin{displaymath}
|\omega|=\frac 1{|\lambda|}\leq \frac 1{C_1}\sqrt[3]{|a|+|c|}\mu^{-\frac 23}.
\end{displaymath}
The regularity of $(\phi,\psi)$ can be obtained by a standard bootstrapping argument (see for instance \cite{Chen2}, Proposition 3.2).

\medskip
\noindent
In order to prove the exponential decay of $\phi$ and $\psi$, following the ideas of  Theorem 8.1.1 in \cite{Cazenave}, we consider, for $\epsilon,\eta>0$, $\displaystyle h(x)=e^{\frac{\epsilon|x|}{1+\eta|x|}}\in L^{\infty}(\er)$. 
Multiplying equations in (\ref{estacionaria}) by $h\phi$ and $h\psi$ respectively and integrating, we get
$$a\int h\phi\phi''+c\int h\psi\psi''+\int h(\phi^2+\psi^2) +\frac 32\int h\phi^2\psi-2\omega\int h\phi\psi=0$$
Integrating by parts and using the fact that $h'\leq \epsilon h$, we obtain
$$\int h(\phi^2+\psi^2-2\omega \phi\psi)\leq a\int h\phi'^2+c\int h\psi'^2+\epsilon\int h(|\phi\phi'|+|\psi\psi'|)+\frac 32\int h\phi^2|\psi|,$$
and
\begin{multline*}
\int h\left(\left(1-\frac{\epsilon}2\right)\phi^2+\left(1-\frac{\epsilon}2\right)\psi^2-2\omega \phi\psi\right)\leq \\
\left(a+\frac{\epsilon}2\right)\int h\phi'^2+\left(c+\frac{\epsilon}2\right)\int h\psi'^2+\frac 32\int h\phi^2|\psi|\leq \frac 32\int h\phi^2|\psi|
\end{multline*}
for $\epsilon$ small enough. 
Since $\psi\in H^1(\er)$, $\displaystyle\lim_{|x|\to +\infty} \psi(x)=0$. For $\epsilon'>0$ to be chosen later, we set $r>0$ such that $|\psi(x)|\leq \epsilon'$ for $|x|>r$.
We then get
$$\int h\left(\left(1-\frac{\epsilon}2-\frac{3\epsilon_1}2\right)\phi^2+\left(1-\frac{\epsilon}2\right)\psi^2-2\omega \phi\psi\right)\leq \frac 32\int_{|x|\leq r}h\phi^2|\psi|$$
and
$$\int h(C_1\phi^2+C_2\psi^2)\leq \frac 32\int_{|x|\leq r}h\phi^2|\psi|$$
where $\displaystyle C_1=1-\frac{\epsilon}2-\frac{3\epsilon_1}2-|\omega|>0$ and 
$\displaystyle C_2=1-\frac{\epsilon}2-|\omega|>0$ for $\epsilon,\epsilon'$ small enough (and for $|\omega|$ small).\\
Finally, taking $\eta\to 0$, by Fatou's Lemma and Lebesgue's Theorem, we obtain
$$\int e^{\epsilon |x|}\phi^2<+\infty\,\textrm{ and }\,\int e^{\epsilon |x|}\phi^2<+\infty.$$
In view of Theorem 8.1.7 of \cite{Cazenave}, this is enough to conclude that 
\begin{displaymath}
e^{\alpha |x|}\phi,e^{\alpha |x|}\psi\in L^{\infty}(\er)
\end{displaymath}
for some $0<\alpha\leq \epsilon$.

\end{document}